\documentclass{amsart}
\usepackage{amssymb,latexsym}
\numberwithin{equation}{section}
\newtheorem{theorem}{Theorem}[section]

\newtheorem{corollary}[theorem]{Corollary}

\newtheorem{example}[theorem]{Example}

\begin{document}
\pagenumbering{arabic}
\pagestyle{headings}
\def\sof{\hfill\rule{2mm}{2mm}}
\def\ls{\leq}
\def\gs{\geq}
\def\SS{\mathcal S}
\def\qq{{\bold q}}
\def\txx{{\frac1{2\sqrt{x}}}}

\title{{\bf coloured permutations containing and avoiding certain patterns}}

\author{Toufik Mansour}
\maketitle
\begin{center}{LABRI, Universit\'e Bordeaux I,\\
              351 cours de la Lib\'eration, 33405 Talence Cedex, France\\
        {\tt toufik@labri.fr} }
\end{center}
%
%===========================================================================
\section*{\bf Abstract}
Following \cite{M2}, let $S_n^{(r)}$ be the set of all coloured
permutations on the symbols $1,2,\dots,n$ with colours
$1,2,\dots,r$, which is the analogous of the symmetric group when
$r=1$, and the hyperoctahedral group when $r=2$. Let
$I\subseteq\{1,2,\dots,r\}$ be subset of $d$ colours; we define
$T_{k,r}^m(I)$ be the set of all coloured permutations $\phi\in
S_k^{(r)}$ such that $\phi_1=m^{(c)}$ where $c\in I$. We prove
that, the number $T_{k,r}^m(I)$-avoiding coloured permutations in
$S_n^{(r)}$ equals $(k-1)!r^{k-1}\prod_{j=k}^n h_j$ for $n\geq k$
where $h_j=(r-d)j+(k-1)d$. We then prove that for any $\phi\in
T_{k,r}^1(I)$ (or any $\phi\in T_{k,r}^k(I)$), the number of
coloured permutations in $S_n^{(r)}$ which avoid all patterns in
$T_{k,r}^1(I)$ (or in $T_{k,r}^k(I)$) except for $\phi$ and
contain $\phi$ exactly once equals $\prod_{j=k}^n h_j\cdot
\sum_{j=k}^n \frac{1}{h_j}$ for $n\geq k$. Finally, for any
$\phi\in T_{k,r}^m(I)$, $2\leq m\leq k-1$, this number equals
$\prod_{j=k+1}^n h_j$ for $n\geq k+1$. These results generalize
recent results due to Mansour \cite{MM}, and due to Simion
\cite{Rs}.
%===================================================================================
\section{\bf Introduction}
The main goal of this note is to give analogies of enumerative
results on certain classes of permutations characterized by
pattern-avoidance in the symmetric group (see \cite{MM}), and in
the hyperoctahedral group (see \cite{Rs}). In $S_n^{(r)}$ (see
\cite{M2}), the natural analogue of the symmetric group and of the
hyperoctahedral group, we identify classes of restricted coloured
permutations with enumerative properties analogous to results in
the symmetric group and hyperoctahedral group. In the remainder of
this section we present a brief account of earlier work which
motivated out investigation, summarize the main results, and
present
the basic definitions used throughout the note.\\

Pattern avoidance in the symmetric group proved to be a useful
language in a variety of seemingly unrelated problems, from stack
sorting \cite{Kn,Rt,W} to the theory of Kazhdan-Lusztig
polynomials ~\cite{Fb}, singularities of Schubert varieties
\cite{LS,SCb}, Chebyshev polynomials \cite{CW,MV1,Kr,MV2,MV3}, and
rook polynomials \cite{MV4}. Signed pattern avoidance in the
hyperoctahedral group proved to be a useful language in
combinatorial statistics defined in type-$B$ noncrossing
partitions, enumerative combinatorics \cite{Rs,MbRs}, algebraic
combinatorics \cite{FK,BK,Be,Cm,Vr}. \\

Let $\pi\in S_n$ and $\tau\in S_k$ be two permutations. An {\it
occurrence} of $\tau$ in $\pi$ is a subsequence $1\leq
i_1<i_2<\dots<i_k\leq n$ such that $(\pi_{i_1},\dots,\pi_{i_k})$
is order-isomorphic to $\tau$; in such a context $\tau$ is usually
called a {\it pattern}. We say that $\pi$ {\it avoids} $\tau$, or
is $\tau$-{\it avoiding}, if there is no occurrence of $\tau$ in
$\pi$. The set of all $\tau$-avoiding permutations in $S_n$ is
denoted by $S_n(\tau)$. For an arbitrary finite collection of
patterns $T$, we say that $\pi$ avoids $T$ if $\pi$ avoids any
$\tau\in T$; the corresponding subset of $S_n$ is denoted by
$S_n(T)$. The first case examined was the case of permutations
avoiding one pattern of length $3$. Knuth \cite{Kn} found that
$|S_n(\tau)|=C_n$ for all $\tau\in S_3$, where $C_n$ is the $n$th
Catalan number. Later, Simion and Schmidt \cite{SS}
found the cardinalities of $|S_n(T)|$ for all $T\subset S_3$. \\

The hyperoctahedral group $B_n$ is an analog of the symmetric
group $S_n$. Let us view the elements of $B_n$ as signed
permutation $b=b_1b_2\dots b_n$ in which each of the symbols
$1,2,\dots,n$ appears once, possibly barred. Thus, the cardinality
of $B_n$ is $n!2^n$. Simion \cite{Rs} was looking for the analogs
of Knuth's results for $B_n$; she discovered that for every
$2$-letter signed pattern $\tau$; the number of $\tau$-avoiding
signed permutations in $B_n$ is $\sum_{j=0}^n {n\choose j}^2j!$.
Besides, Simion \cite{Rs} found the number of all signed
permutations in $B_n$ avoiding double $2$-letter signed patterns
in $B_2$. This invites us to define a further generalizations for
avoiding a pattern in $S_n$ and avoiding a signed pattern in $B_n$.\\

Following \cite{M2} (see also \cite{St}), the group
$S_n^{(r)}=S_n\wr C_r$ where $C_r$ is the cyclic group of order
$r$, is an analog of the symmetric group ($S_n$) and of the
hyperoctahedral group ($B_n$). We will view the elements of the
set $S_n^{(r)}$ as coloured permutations
$\phi=(\phi_1,\phi_2,\dots,\phi_n)$ in which each of the symbols
$1,2,\dots,n$ appears once, coloured by one of the colours
$1,2,\dots,r$ (more generally, we denote by
$S_{\{a_1,\dots,a_n\}}^{\{s_1,\dots,s_r\}}$ the set of all
permutations of the symbols $a_1,\dots,a_n$ where each symbol
appears once and is coloured by one of the colours
$s_1,\dots,s_r$). Thus, $S_n^{(1)}=S_n$, $S_n^{(2)}=B_n$, and the
cardinality of $S_n^{(r)}$ is $n!r^n$. The absolute value notation
means $|\phi|$ is the permutation $(|\phi_1|,\dots,|\phi_n|)$
where $|\phi_j|$ is the symbol which appear in $\phi$ at the
position $j$. An example $\phi=(1^{(1)},3^{(2)},2^{(1)})$ is a
coloured permutation in $S_3^{(2)}$, and $|\phi|=(1,3,2)$.\\

Let $\phi=(\tau_1^{(s_1)},\dots,\tau_k^{(s_k)})\in S_k^{(r)}$, and
$\psi=(\alpha_1^{(v_1)},\dots,\alpha_n^{(v_n)})\in S_n^{(r)}$; we
say that $\psi$ {\em contains} $\phi$ (or is $\phi$-containing) if
there is a sequence of $k$ indices, $1\leq i_1<i_2<\dots<i_k\leq
n$ such that the following two conditions hold:
\begin{enumerate}
\item[(i)]  $(\alpha_{i_1},\dots,\alpha_{i_k})$ is order-isomorphic to
$|\phi|$;\\

\item[(ii)]     $v_{i_j}=s_j$ for all $j=1,2,\dots,k$.
\end{enumerate}

Otherwise, we say that $\psi$ {\em avoids} $\phi$ (or is
$\phi$-avoiding). The set of all $\phi$-avoiding coloured
permutations in $S_n^{(r)}$ is denoted by $S_n^{(r)}(\phi)$, and
in this context $\phi$ is called a {\it coloured pattern\/}. For
an arbitrary finite collection of coloured patterns $T$, we say
that $\psi$ avoids $T$ if $\psi$ avoids any $\phi\in T$; the
corresponding subset of $S_n^{(r)}$ is denoted by $S_n^{(r)}(T)$.
As an example, $\psi=(1^{(1)},2^{(2)},3^{(2)})\in S_3^{(2)}$
avoids $(1^{(1)},2^{(1)})$; that is, $\psi\in S_3^{(2)}((1^{(1)},2^{(1)}))$.\\

In this note, we present an analogs of Mansour's results for
avoiding and containing certain patterns in $S_n$ (see \cite{MM}),
and Simion's results for avoiding signed patterns in $B_n$ (see
\cite[Sec. 3]{Rs}).
%===================================================================================
\section{\bf Coloured permutations avoiding $T_k^m(I_d)$}
Let $I_d$ be any subset of $\{1,2,\dots,r\}$ of $d$ elements, and
let us define $T_{k,r}^m(I_d)$ be the set of all coloured
permutations $\phi\in S_k^{(r)}$ such that $\phi_1=m^{(c)}$ where
$c\in I_d$; that is,
    $$T_{k,r}^m(I_d)=\bigcup\limits_{c\in I_d} \{ \phi\in S_k^{(r)} | \phi_1=m^{(c)}\}.$$

\begin{theorem}
\label{thm1}
Let $k,r\geq 1$, $n\geq k$, and $k\geq m\geq 1$. Then
    $$|S_n^{(r)}(T_{k,r}^m(I_d))|=(k-1)!r^{k-1}\prod_{j=k}^n ( (r-d)j+(k-1)d).$$
\end{theorem}
\begin{proof}
  Let $G_n=S_n^{(r)}(T_{k,r}^m(I_d))$, and define the family of functions
$f_{h,c}:S_n^{(r)}\rightarrow S_{n+1}^{(r)}$ by:
  $$[f_{h,c}(\phi)]_i=\left\{ \begin{array}{ll} h^{(c)},& \mbox{when}\ i=1 \\
                       \phi_{i-1},  & \mbox{when}\ |\phi_{i-1}|<h \\
                       (|\phi_{i-1}|+1)^{(a_{i-1})},& \mbox{when}\ |\phi_{i-1}|\geq h \end{array} \right.$$
for every $i=1,...,n+1$, $\phi\in S_n^{(r)}$, $1\leq c\leq r$ and
$h=1,...,n+1$, where $a_i$ is the colour of the symbol $|\phi_i|$
in $\phi$.\\

From this we see that if $\phi\in G_n$, then
     $$f_{n+1,c_j}(\phi ),f_{n,c_j}(\phi ),...,f_{n+m-k+2,c_j}(\phi),f_{1,c_j}(\phi),...,f_{m-1,c_j}(\phi)\in G_{n+1}$$
for all $j=1,2,\dots,d$, and $f_{h,c}(\phi)\in G_{n+1}$ for all $c\notin I_d$ and
$h=1,2,\dots,n+1$.
so $((k-1)d+(r-d)(n+1)) \cdot |G_n|\leq |G_{n+1}|$ where $n\geq k$.\\

Assume that $((k-1)d+(r-d)(n+1))\cdot |G_n|<|G_{n+1}|$. Then there
exists a coloured permutation $\psi\in G_{n+1}$ such that $m\leq
|\psi_1|\leq n+m-k+1$ and the symbol $|\psi_1|$ coloured by $c\in
I_d$, so there exist $k-1$ positions $1<i_1<\cdots <i_{k-1}\leq
n+1$ such that the subsequence
$\phi_1,\phi_{i_1},\dots,\phi_{i_{k-1}}$ is contains one of the
patterns in $T_{k,r}^m(I_d)$, which contradicts the definition of
$G_{n+1}$. So $((k-1)d+(r-d)(n+1))\cdot |G_n|=|G_{n+1}|$ for
$n\geq k$. Besides $|G_k|=(rk-d)(k-1)!r^{k-1}$ (from the
definitions), hence the theorem holds.
\end{proof}

\begin{example} {\rm(see \cite[Th. 1]{MM})}
Let $T_k^m=T_{k,1}^m(1)$ and $k\geq m\geq 1$. Theorem \ref{thm1} yields
for all $n\geq k$ that
        $$|S_n(T_k^m)|=(k-1)!(k-1)^{n-k+1}.$$
\end{example}

\begin{example} {\rm( see \cite[Sec. 3]{Rs})}
Theorem \ref{thm1} yields for all $n\geq k\geq m\geq 1$ that
        $$|S_n^{(2)}(T_{k,2}^m(1))|=\frac{(n+k-1)!}{\prod_{i=1}^{k-1} (2i-1)}.$$
So, for $k=2$ we have \cite[Eq. $47$]{Rs}.
\end{example}

\begin{example}{\rm(see \cite[Cor. $4.2$]{M2})}
Theorem \ref{thm1} yields for $k=2$, $m=1$, and $I_d=\{1\}$ that
    $$|S_n^{(r)}(T_{2,r}^1(1))|=\prod_{j=0}^n(1+j(r-1)).$$
\end{example}

\begin{corollary}
\label{cm1}
Let $k,r\geq 1$, and $k\geq b\geq a\geq 1$. For $n\geq k$
     $$|S_n^{(r)}(\cup_{m=a}^b T_{k,r}^m(I_d) )|=(k-1)!r^{k-1}\prod_{j=k}^n (d(k+a-b-1)+j(r-d)).$$
\end{corollary}
\begin{proof}
Let $G_n=S_n^{(r)}(T_{k,r}^m(I_d))$.
From Theorem \ref{thm1} we get that $\phi\in G_n$ if and only if either
  $$f_{1,c_j}(\phi),\dots,f_{a-1,c_j}(\phi),f_{n+b-(k-2),c_j}(\phi),\dots,f_{n+1,c_j}(\phi) \in G_{n+1}$$
for $j=1,2,\dots,d$, or $f_{h,c}(\phi)\in G_{n+1}$ for $h=1,2,\dots,n+1$,
$c\notin I_d$.
So $|G_{n+1}|=(d(k+a-b-1)+(r-d)(n+1))|G_n|$. Besides
$|G_k|=(rk+d(a-b-1))(k-1)!r^{k-1}$, hence the theorem holds.
\end{proof}

\begin{example}
Corollary \ref{cm1} yields for all $n\geq 0$
    $$|S_n^{(2)}( T_{k,2}^1(1)\cup T_{k,2}^2(1) )|=2^{k-1}n!.$$
\end{example}

This example invite us to generalize. By use Corollary \ref{cm1} we get the following.

\begin{corollary}
\label{cm1}
Let $k,r\geq 1$, and $n\geq k$. Then
     $$|S_n^{(r)}(\cup_{m=1}^k T_{k,r}^m(I_d) )|=r^{k-1}(r-d)^{n+1-k}n!.$$
\end{corollary}
%=============================================================================
\section{\bf Avoiding $T_{k,r}^1(I_d)\backslash{\{\phi\}}$ and containing $\phi$ exactly once}

Let $M_{k,r}^{m}(\phi,I_d)=T_{k,r}^m(I_d\backslash\{\phi\}$,
for $\phi\in T_{k,r}^m(I_d)$. We denote by $S_n^{(r)}(T_{k,r}^m(I_d);\tau)$
the set of all permutations in $S_n^{(r)}$ that avoid $M_{k,r}^m(\phi,I_d)$
and contain $\phi$ exactly once.

\begin{theorem}
\label{thm2}
Let $k,r\geq 1$, and $n\geq k$. Then
    $$|S_n^{(r)}(T_{k,r}^1(I_d);\phi)|=\prod_{j=k}^n (d(k-1)+(r-d)j)\cdot \sum_{j=k}^n \frac{1}{d(k-1)+(r-d)j},$$
for all $\phi\in T_{k,r}^1(I_d)$.
\end{theorem}
\begin{proof}
Let $\alpha\in S_n^{(r)}(T_{k,r}^1(I_d);\phi)$, and let us consider the
possible values of $\psi_1$:
\begin{enumerate}
  \item $|\psi_1|$ is coloured by colour $c\notin I_d$.
    Evidently $\psi\in S_n^{(r)}(T_{k,r}^1(I_d);\phi)$
    if and only if $\phi$ avoids $M_{k,r}^m(\phi,I_d)$ and contains $\phi$
    exactly once.\\

  \item $|\psi_1|\geq n-k+2$ is coloured by $c\in I_d$.
    Evidently $\psi\in S_n^{(r)}(T_{k,r}^1(I_d);\phi)$
    if and only if $\phi$ avoids $M_{k,r}^m(\phi,I_d)$ and contains $\phi$
    exactly once.\\

  \item $|\psi_1|\leq n-k$ is coloured by $c\in I_d$.
    Then there exist $1<i_1<\dots<i_k\leq n$ such
        that $(\psi_1,\psi_{i_1},\dots,\psi_{i_k})$ is a
    coloured permutation of the symbols $n,\dots,n-k+1,\alpha_1$ coloured by any colours
    such that $\psi_1$ coloured by $c\in I_d$.
    For any choice of $k-1$ positions out of $i_1,\dots,i_k$, the
    corresponding coloured permutations proceeded by $|\psi_1|$, is
    containing some coloured pattern in $T_{k,r}^1(I_d)$. Since
    $\psi$ avoids $M_{k,r}^m(\phi,I_d)$, it is, in fact,
    order-isomorphic to $\phi$. We thus get at least $k$ occurrences
    of $\phi$ in $\psi$, a contradiction.\\

  \item $|\psi_1|=n-k+1$. Then there exist $1<i_1<\dots<i_{k-1}\leq n$ such that
    $\eta =(\psi_1,\psi_{i_1},\dots,\psi_{i_{k-1}})$ is a
    coloured permutation of the symbols $n,\dots,n-k+1$. As above, we
    immediately get that $\eta$ is order-isomorphic to $\phi$.
    We denote by $A_n$ the set of all coloured permutations
    in $S_n^{(r)}(T_{k,r}^1(I_d);\phi)$
    such that the symbol $|\phi_1|$ is $n-k+1$ and coloured by $c\in I_d$,
    and define the family of functions $f_h:A_n\rightarrow S_{n+1}^{(r)}$ by:
        $$[f_{h,c}(\beta)]_i=\left\{ \begin{array}{ll} 1^{(c)},  & \mbox{when}\ i=h \\
                       (|\beta_i|+1)^{(a_i)},     & \mbox{when}\ i<h \\
                       (|\beta_{i-1}|+1)^{(a_i)}, & \mbox{when}\ i>h \end{array} \right. ,$$
    for every $i=1,...,n+1$, $\beta\in A_n$, $h=1,...,n+1$ and $c=1,\dots,r$, where
    the symbol $|\beta_i|$ coloured by the colour $a_i$.
    It is easy to see that for all $\beta\in A_n$,
         $$f_{n+1,c_j}(\beta),\dots,f_{n-k+3,c_j}(\beta)\in A_{n+1},$$
    and $f_{h,c}(\beta)\in A_{n+1}$ for $c\notin I_d$ and $h=2,\dots,n$
    (for $h=1$ we added in the first case),
    hence $(d(k-1)+(n)(r-d))|A_n|\leq |A_{n+1}|$.
        Now we define another function $g:A_{n+1}\rightarrow S_n^{(r)}$ by:
        $$[g(\beta)]_i=\left\{ \begin{array}{ll}
                (|\beta_i|-1)^{(a_i)},     & \mbox{when}\ i<h \\
                            (|\beta_{i+1}|-1)^{(a_i)}, & \mbox{when}\ i+1>h \end{array} \right. ,$$
    where $|\beta_h|=1$, $i=1,...,n$, $\beta\in A_{n+1}$.\\

    Observe that $h\geq n-k+3$ such that $|\beta_h|$ coloured by $c\in I_d$,
    since otherwise already
    $(\beta_h,\beta_{h+1},\dots,\beta_{n+1})$ contains a pattern
    from $T_k^1(I_d)$, a contradiction. It is easy to see that
    $g(\beta)\in A_n$ for all $\beta\in A_{n+1}$, hence
    $|A_{n+1}|\leq (d(k-1)+(r-d)(n+1))|A_n|$. So finally,
    $|A_{n+1}|=(d(k-1)+(r-d)n)|A_n|$ and
    $|A_n|=\prod_{j=k}^{n-1} (d(k-1)+(r-d)j)$, since $|A_k|=1$.
   \end{enumerate}
Since the above cases $1,2,3,4$ are disjoint we obtain
    $$
    |S_n^{(r)}(T_{k,r}^1(I_d);\phi)|=(d(k-1)+(r-d)n)|S_{n-1}^{(r)}(T_{k,r}^1(I_d);\phi)|+\prod_{j=k}^{n-1} (d(k-1)+(r-d)j).
    $$
Besides $|S_k^{(r)}(T_{k,r}^1(I_d);\phi)|=1$, hence the theorem holds.
\end{proof}

\begin{example}
By Theorem \ref{thm2} we get as follows. For any $\phi\in T_{k,2}^1(1)$
    $$|S_n^{(2)}(T_{k,2}^1(1);\phi )|=\frac{(n+1)!}{(2k-2)!}\sum_{j=2k-1}^{n+k-1} \frac{1}{j}.$$
\end{example}

By the natural bijection between the set
$S_n^{(r)}(T_{k,r}^1(I_d);\phi)$ and the set
$S_n^{(r)}(T_{k,r}^k(I_d);\phi')$ for all $\phi\in
T_{k,r}^1(I_d)$, where $\phi'_i=(k+1-|\phi_i|)^{(a_i)}$, $a_i$ is
colour of the symbol $|\phi_i|$ in $\phi$, $i=1,2,\dots,k$, hence
by Theorem \ref{thm2} we get the following.

\begin{corollary}
\label{cm2}
Let $k,r\geq 1$, and $n\geq k$. Then
   $$|S_n^{(r)}(T_{k,r}^k(I_d);\phi)|=\prod_{j=k}^n (d(k-1)+(r-d)j)\cdot \sum_{j=k}^n \frac{1}{d(k-1)+(r-d)j},$$
for all $\phi\in T_{k,r}^k(I_d)$.
\end{corollary}
%=============================================================================
\section{\bf Avoiding $T_{k,r}^m(I_d)\backslash{\{\phi\}}$ and containing $\phi$ exactly once, $2\leq m \leq k-1$}

Now we calculate the cardinalities of the sets
$S_n^{(r)}(T_{k,r}^m(I_d);\phi)$ where $2\leq m\leq k-1$, $\phi\in T_{k,r}^m(I_d)$.

\begin{theorem}
\label{thm3}
Let $r\geq 1$, $k\geq 3$, and $n\geq k+1$. Then
   $$|S_n^{(r)}(T_{k,r}^m(I_d);\phi)|=\prod_{j=k+1}^n (d(k-1)+(r-d)j),$$
   for all $2\leq m\leq k-1$, $\phi\in T_{k,r}^m(I_d)$.
\end{theorem}
\begin{proof}
Let $G_n=S_n^{(r)}(T_{k,r}^m(I_d);\phi)$, $\psi\in G_n$, and let us consider the
possible values of $\psi_1$:
\begin{enumerate}
  \item $|\psi_1|\leq m-1$ and coloured by $c\in I_d$.
    Evidently $\psi\in G_n$ if and only if
    $(\psi_2,\dots,\psi_n)$ avoids $M_{k,r}^m(\phi,I_d)$ and contains
    $\phi$ exactly once.\\

  \item $|\psi_1\geq n-k+m+1$ and coloured by $c\in I_d$.
    Evidently $\psi\in G_n$ if and only if
    $(\psi_2,\dots,\psi_n)$ avoids $M_{k,r}^m(\phi,I_d)$ and contains
    $\phi$ exactly once.\\

  \item $m\leq|\psi_1|\leq n-k+m$ and coloured by one of the colours $c_1,\dots,c_d$.
    By definition we have that $|G_k|=1$, so
    let $n\geq k+1$. If $m+1\leq|\psi_1|$ then $\psi$ contains at
    least $m\geq 2$ occurrences of a pattern from $T_{k,r}^m(I_d)$, and
    If $|\psi_1|\leq n-k+m-1$ then $\psi$ contains at least
    $k-m+1\geq 2$ occurrences of a pattern from $T_{k,r}^m(I_d)$,
    a contradiction.

  \item $|\psi_1|$ coloured by a colour $c\notin I_d$.
    Evidently $\psi\in G_n$ if and only if
    $(\psi_2,\dots,\psi_n)$ avoids $M_{k,r}^m(\phi, I_d)$ and
    contains $\phi$ exactly once.\\
  \end{enumerate}
Since the above cases $1,2,3,4$ are disjoint we obtain $|G_n|=(d(k-1)+(r-d)n)|G_{n-1}|$ for all $n\geq k+1$.
Besides $|G_k|=1$, hence the theorem holds.
\end{proof}

\begin{example}{\rm (see \cite{Rs})}
By Theorem \ref{thm3} we get as follows. For any $\phi\in T_{k,2}^m(1)$
    $$|S_n^{(2)}(T_{k,2}^m(1);\phi)|=\frac{(n+k-1)!}{(2k-1)!},$$
where $2\leq m\leq k-1$, $k\geq 3$ and $n\geq k$.
\end{example}
%=============================================================================

\end{document}